\documentclass[reqno,12pt,a4paper]{amsart}

\usepackage{mathrsfs}
\usepackage{amssymb}
\usepackage{amsfonts}
\usepackage{latexsym}
\usepackage{amsthm}
\usepackage{graphicx}
\def\lb{\label}

\newcommand{\er}[1]{\textrm{(\ref{#1})}}

\begin{document}

%%%%%%%%%% Some definitions %%%%%%%%%%

%%%%%%%% Equations, theorems %%%%%%%%%
\renewcommand{\theequation}{\arabic{section}.\arabic{equation}}
\theoremstyle{plain}
\newtheorem{theorem}{\bf Theorem}[section]
\newtheorem{lemma}[theorem]{\bf Lemma}
\newtheorem{corollary}[theorem]{\bf Corollary}
\newtheorem{proposition}[theorem]{\bf Proposition}
\newtheorem{definition}[theorem]{\bf Definition}
\newtheorem{remark}[theorem]{\bf Remark}

%%%%% Alphabet %%%%%
\def\a{\alpha}  \def\cA{{\mathcal A}}     \def\bA{{\bf A}}  \def\mA{{\mathscr A}}
\def\b{\beta}   \def\cB{{\mathcal B}}     \def\bB{{\bf B}}  \def\mB{{\mathscr B}}
\def\g{\gamma}  \def\cC{{\mathcal C}}     \def\bC{{\bf C}}  \def\mC{{\mathscr C}}
\def\G{\Gamma}  \def\cD{{\mathcal D}}     \def\bD{{\bf D}}  \def\mD{{\mathscr D}}
\def\d{\delta}  \def\cE{{\mathcal E}}     \def\bE{{\bf E}}  \def\mE{{\mathscr E}}
\def\D{\Delta}  \def\cF{{\mathcal F}}     \def\bF{{\bf F}}  \def\mF{{\mathscr F}}
\def\c{\chi}    \def\cG{{\mathcal G}}     \def\bG{{\bf G}}  \def\mG{{\mathscr G}}
\def\z{\zeta}   \def\cH{{\mathcal H}}     \def\bH{{\bf H}}  \def\mH{{\mathscr H}}
\def\e{\eta}    \def\cI{{\mathcal I}}     \def\bI{{\bf I}}  \def\mI{{\mathscr I}}
\def\p{\psi}    \def\cJ{{\mathcal J}}     \def\bJ{{\bf J}}  \def\mJ{{\mathscr J}}
\def\vT{\Theta} \def\cK{{\mathcal K}}     \def\bK{{\bf K}}  \def\mK{{\mathscr K}}
\def\k{\kappa}  \def\cL{{\mathcal L}}     \def\bL{{\bf L}}  \def\mL{{\mathscr L}}
\def\l{\lambda} \def\cM{{\mathcal M}}     \def\bM{{\bf M}}  \def\mM{{\mathscr M}}
\def\L{\Lambda} \def\cN{{\mathcal N}}     \def\bN{{\bf N}}  \def\mN{{\mathscr N}}
\def\m{\mu}     \def\cO{{\mathcal O}}     \def\bO{{\bf O}}  \def\mO{{\mathscr O}}
\def\n{\nu}     \def\cP{{\mathcal P}}     \def\bP{{\bf P}}  \def\mP{{\mathscr P}}
\def\r{\rho}    \def\cQ{{\mathcal Q}}     \def\bQ{{\bf Q}}  \def\mQ{{\mathscr Q}}
\def\s{\sigma}  \def\cR{{\mathcal R}}     \def\bR{{\bf R}}  \def\mR{{\mathscr R}}
\def\S{\Sigma}  \def\cS{{\mathcal S}}     \def\bS{{\bf S}}  \def\mS{{\mathscr S}}
\def\t{\tau}    \def\cT{{\mathcal T}}     \def\bT{{\bf T}}  \def\mT{{\mathscr T}}
\def\f{\phi}    \def\cU{{\mathcal U}}     \def\bU{{\bf U}}  \def\mU{{\mathscr U}}
\def\F{\Phi}    \def\cV{{\mathcal V}}     \def\bV{{\bf V}}  \def\mV{{\mathscr V}}
\def\P{\Psi}    \def\cW{{\mathcal W}}     \def\bW{{\bf W}}  \def\mW{{\mathscr W}}
\def\o{\omega}  \def\cX{{\mathcal X}}     \def\bX{{\bf X}}  \def\mX{{\mathscr X}}
\def\x{\xi}     \def\cY{{\mathcal Y}}     \def\bY{{\bf Y}}  \def\mY{{\mathscr Y}}
\def\X{\Xi}     \def\cZ{{\mathcal Z}}     \def\bZ{{\bf Z}}  \def\mZ{{\mathscr Z}}
\def\O{\Omega}

%*********************
\def\be{{\bf e}}
\def\bv{{\bf v}} \def\bu{{\bf u}}
 \def\mn{\mathrm n}
\def\mm{\mathrm m}
%************************

\newcommand{\mc}{\mathscr {c}}

\newcommand{\gA}{\mathfrak{A}}          \newcommand{\ga}{\mathfrak{a}}
\newcommand{\gB}{\mathfrak{B}}          \newcommand{\gb}{\mathfrak{b}}
\newcommand{\gC}{\mathfrak{C}}          \newcommand{\gc}{\mathfrak{c}}
\newcommand{\gD}{\mathfrak{D}}          \newcommand{\gd}{\mathfrak{d}}
\newcommand{\gE}{\mathfrak{E}}
\newcommand{\gF}{\mathfrak{F}}           \newcommand{\gf}{\mathfrak{f}}
\newcommand{\gG}{\mathfrak{G}}           %\newcommand{\gg}{\mathfrak{g}}
\newcommand{\gH}{\mathfrak{H}}           \newcommand{\gh}{\mathfrak{h}}
\newcommand{\gI}{\mathfrak{I}}           \newcommand{\gi}{\mathfrak{i}}
\newcommand{\gJ}{\mathfrak{J}}           \newcommand{\gj}{\mathfrak{j}}
\newcommand{\gK}{\mathfrak{K}}            \newcommand{\gk}{\mathfrak{k}}
\newcommand{\gL}{\mathfrak{L}}            \newcommand{\gl}{\mathfrak{l}}
\newcommand{\gM}{\mathfrak{M}}            \newcommand{\gm}{\mathfrak{m}}
\newcommand{\gN}{\mathfrak{N}}            \newcommand{\gn}{\mathfrak{n}}
\newcommand{\gO}{\mathfrak{O}}
\newcommand{\gP}{\mathfrak{P}}             \newcommand{\gp}{\mathfrak{p}}
\newcommand{\gQ}{\mathfrak{Q}}             \newcommand{\gq}{\mathfrak{q}}
\newcommand{\gR}{\mathfrak{R}}             \newcommand{\gr}{\mathfrak{r}}
\newcommand{\gS}{\mathfrak{S}}              \newcommand{\gs}{\mathfrak{s}}
\newcommand{\gT}{\mathfrak{T}}             \newcommand{\gt}{\mathfrak{t}}
\newcommand{\gU}{\mathfrak{U}}             \newcommand{\gu}{\mathfrak{u}}
\newcommand{\gV}{\mathfrak{V}}             \newcommand{\gv}{\mathfrak{v}}
\newcommand{\gW}{\mathfrak{W}}             \newcommand{\gw}{\mathfrak{w}}
\newcommand{\gX}{\mathfrak{X}}               \newcommand{\gx}{\mathfrak{x}}
\newcommand{\gY}{\mathfrak{Y}}              \newcommand{\gy}{\mathfrak{y}}
\newcommand{\gZ}{\mathfrak{Z}}             \newcommand{\gz}{\mathfrak{z}}

\def\ve{\varepsilon}   \def\vt{\vartheta}    \def\vp{\varphi}    \def\vk{\varkappa}

\def\A{{\mathbb A}} \def\B{{\mathbb B}} \def\C{{\mathbb C}}
\def\dD{{\mathbb D}} \def\E{{\mathbb E}} \def\dF{{\mathbb F}} \def\dG{{\mathbb G}} \def\H{{\mathbb H}}\def\I{{\mathbb I}} \def\J{{\mathbb J}} \def\K{{\mathbb K}} \def\dL{{\mathbb L}}\def\M{{\mathbb M}} \def\N{{\mathbb N}} \def\dO{{\mathbb O}} \def\dP{{\mathbb P}} \def\R{{\mathbb R}}\def\S{{\mathbb S}} \def\T{{\mathbb T}} \def\U{{\mathbb U}} \def\V{{\mathbb V}}\def\W{{\mathbb W}} \def\X{{\mathbb X}} \def\Y{{\mathbb Y}} \def\Z{{\mathbb Z}}

%%%%% Arrows %%%%%

\def\la{\leftarrow}              \def\ra{\rightarrow}            \def\Ra{\Rightarrow}
\def\ua{\uparrow}                \def\da{\downarrow}
\def\lra{\leftrightarrow}        \def\Lra{\Leftrightarrow}

%%%%% Typography %%%%%

\def\lt{\biggl}                  \def\rt{\biggr}
\def\ol{\overline}               \def\wt{\widetilde}
\def\no{\noindent}

%%%%% Math signs %%%%%

\let\ge\geqslant                 \let\le\leqslant
\def\lan{\langle}                \def\ran{\rangle}
\def\/{\over}                    \def\iy{\infty}
\def\sm{\setminus}               \def\es{\emptyset}
\def\ss{\subset}                 \def\ts{\times}
\def\pa{\partial}                \def\os{\oplus}
\def\om{\ominus}                 \def\ev{\equiv}
\def\iint{\int\!\!\!\int}        \def\iintt{\mathop{\int\!\!\int\!\!\dots\!\!\int}\limits}
\def\el2{\ell^{\,2}}             \def\1{1\!\!1}
\def\sh{\sharp}
\def\wh{\widehat}
\def\bs{\backslash}
\def\intl{\int\limits}
%%%%% Math operations %%%%%

\def\na{\mathop{\mathrm{\nabla}}\nolimits}
\def\sh{\mathop{\mathrm{sh}}\nolimits}
\def\ch{\mathop{\mathrm{ch}}\nolimits}
\def\where{\mathop{\mathrm{where}}\nolimits}
\def\all{\mathop{\mathrm{all}}\nolimits}
\def\as{\mathop{\mathrm{as}}\nolimits}
\def\Area{\mathop{\mathrm{Area}}\nolimits}
\def\arg{\mathop{\mathrm{arg}}\nolimits}
\def\const{\mathop{\mathrm{const}}\nolimits}
\def\det{\mathop{\mathrm{det}}\nolimits}
\def\diag{\mathop{\mathrm{diag}}\nolimits}
\def\diam{\mathop{\mathrm{diam}}\nolimits}
\def\dim{\mathop{\mathrm{dim}}\nolimits}
\def\dist{\mathop{\mathrm{dist}}\nolimits}
\def\Im{\mathop{\mathrm{Im}}\nolimits}
\def\Iso{\mathop{\mathrm{Iso}}\nolimits}
\def\Ker{\mathop{\mathrm{Ker}}\nolimits}
\def\Lip{\mathop{\mathrm{Lip}}\nolimits}
\def\rank{\mathop{\mathrm{rank}}\limits}
\def\Ran{\mathop{\mathrm{Ran}}\nolimits}
\def\Re{\mathop{\mathrm{Re}}\nolimits}
\def\Res{\mathop{\mathrm{Res}}\nolimits}
\def\res{\mathop{\mathrm{res}}\limits}
\def\sign{\mathop{\mathrm{sign}}\nolimits}
\def\span{\mathop{\mathrm{span}}\nolimits}
\def\supp{\mathop{\mathrm{supp}}\nolimits}
\def\Tr{\mathop{\mathrm{Tr}}\nolimits}
\def\BBox{\hspace{1mm}\vrule height6pt width5.5pt depth0pt \hspace{6pt}}

%%%%%%%%%%%%% specialities %%%%%%%%%%%%%%

\newcommand\nh[2]{\widehat{#1}\vphantom{#1}^{(#2)}}
%{{\mathop{#1}\limits^\wedge}\vphantom{#1}^{(#2)}}
\def\dia{\diamond}

\def\Oplus{\bigoplus\nolimits}

%%%%%%%%%%% End of definitions %%%%%%%%%%

%%%%% OLD OLD OLD

\def\qqq{\qquad}
\def\qq{\quad}
\let\ge\geqslant
\let\le\leqslant
\let\geq\geqslant
\let\leq\leqslant
\newcommand{\ca}{\begin{cases}}
\newcommand{\ac}{\end{cases}}
\newcommand{\ma}{\begin{pmatrix}}
\newcommand{\am}{\end{pmatrix}}
\renewcommand{\[}{\begin{equation}}
\renewcommand{\]}{\end{equation}}
\def\eq{\begin{equation}}
\def\qe{\end{equation}}
\def\[{\begin{equation}}
\def\bu{\bullet}
\def\bq{\mathbf q}

\title[{Estimates of bands for Laplacians on periodic equilateral metric graphs}]
{Estimates of bands for Laplacians on periodic equilateral metric graphs}

\date{\today}

\author[Evgeny Korotyaev]{Evgeny Korotyaev}
\address{Mathematical Physics Department, Faculty of Physics, Ulianovskaya 2,
St. Petersburg State University, St. Petersburg, 198904,
\ korotyaev@gmail.com,}
\author[Natalia Saburova]{Natalia Saburova}
\address{Department of Mathematical Analysis, Algebra and Geometry, Institute of
Mathematics, Information and Space Technologies, Uritskogo St. 68,
Northern (Arctic) Federal University, Arkhangelsk, 163002,
 \ n.saburova@gmail.com}


\begin{abstract}
We consider Laplacians on periodic equilateral metric graphs.
The spectrum of the Laplacian consists of an
absolutely continuous part (which is a union of an infinite number
of non-degenerated spectral bands) plus an infinite number of flat
bands, i.e., eigenvalues of infinite multiplicity. We estimate the Lebesgue measure of the bands on a finite interval in terms of geometric parameters of
the graph.  The proof is based on spectral properties of discrete Laplacians.
\end{abstract}

\maketitle

%%%%%%%%%%%%%%%%%%%%%%%%%%%%%%%%%%%%%%%%%%%%%%%%%%%%%

\vskip 0.25cm

\section {Introduction and main results}
\setcounter{equation}{0}

We consider metric Laplacians $\D_M$ on $\Z^d$-periodic equilateral
metric graphs (each edge has  unit length). Such operators
arise naturally as simplified models in mathematics, physics,
chemistry, and engineering when one considers propagation of waves
of various nature through a quasi-one-dimensional system that looks
like a thin neighborhood of a graph.

It is well-known that the spectrum of the Laplacian $\D_M$ consists of an absolutely continuous part plus an
infinite number of flat bands (i.e., eigenvalues with infinite
multiplicity). These and other properties of $\D_M$ are discussed, e.g., in \cite{BKu12}, \cite{P12} and see references therein. The absolutely continuous spectrum consists of an
infinite number of spectral bands separated by gaps.
There is a known problem: to estimate the spectrum and gaps of the Laplacian on periodic metric graph. Note that in the case of the Schr\"odinger
operators $-\D+Q$ with a periodic potential $Q$ in $\R^d$ there are
two-sided estimates of potentials in terms of gap lengths only at
$d=1$ in \cite{K98}, \cite{K03}. We do not know other estimates.
For the case of periodic
graphs we know only two papers about estimates of spectrum and gaps:

(1) Lled\'o and Post \cite{LP08} estimated the positions of spectral bands of  Laplacians both on metric and discret graphs in terms of eigenvalues of the operator on finite graphs (the so-called eigenvalue bracketing).

(2) Korotyaev and Saburova \cite{KS14} considered Schr\"odinger
operators on the discrete graphs and estimated the Lebesgue measure
of their spectrum  in terms of geometric parameters of the graph
only.

Our main goal is to estimate the spectral bands and gaps for the Laplacian
on the metric graph in terms of geometric parameters of
the graph. Due to Cattaneo \cite{C97}, in order to
study the spectrum of the Laplacian $\D_M$ on  the equilateral metric graph it
is enough to consider only the spectral interval $[0,\pi^2]$.

\subsection{Metric Laplacians.} Let $\G=(V,\cE)$ be a connected graph,
possibly  having loops and multiple edges, where $V$ is the set of
its vertices and  $\cE$ is the set of its unoriented edges. The
graphs under consideration are embedded into $\R^d$. An edge
connecting vertices $u$ and $v$ from $V$ will be denoted as the
unordered pair $(u,v)_e\in\cE$ and is said to be \emph{incident} to
the vertices. Vertices $u,v\in V$ will be called \emph{adjacent} and
denoted by $u\sim v$, if $(u,v)_e\in \cE$. We define the degree
${\vk}_v$ of the vertex $v\in V$ as the number of all its
incident edges from $\cE$ (here a loop is counted twice). Below we consider locally finite
$\Z^d$-periodic metric equilateral graphs $\G$, i.e., graphs satisfying the
following conditions:

1) {\it the number of vertices from $V$ in any bounded domain $\ss\R^d$ is
finite;

2) the degree of each vertex is finite;

3) there exists a basis $a_1,\ldots,a_d$ in $\R^d$ such that $\G$
 is invariant under translations through the vectors $a_1,\ldots,a_d$:
$$
\G+a_s=\G, \qqq  \forall\,s\in\N_d=\{1,\ldots,d\}.
$$
The vectors $a_1,\ldots,a_d$ are called the periods of $\G$.}

4) \emph{All edges of the graph have the unit length.}

In the space $\R^d$ we consider a coordinate system with the origin at
some point $O$. The coordinate axes of this system are directed along
the vectors $a_1,\ldots,a_d$. Below the coordinates of all vertices
of $\G$ will be expressed  in this coordinate system.
From the definition it follows that a $\Z^d$-periodic graph $\G$
is invariant under translations through any integer vector:
$$
\G+m=\G,\qqq \forall\, m\in\Z^d.
$$

Each edge $\be$ of $\G$ will be identified with the segment $[0,1]$.
This identification introduces a local coordinate $t\in[0,1]$ along
each edge. Thus, we give an orientation on the edge. Note that the
spectrum of Laplacians on metric graphs does not depend on the
orientation of graph edges. For each function $y$ on $\G$ we define
a function $y_{\be}=y\big|_{\be}$, $\be\in\cE$. We identify each
function  $y_{\be}$ on $\be$ with a function on $[0,1]$ by using the
local coordinate $t\in[0,1]$.
Let $L^2(\G)$ be the Hilbert space of all function $y=(y_\be)_{\be\in\cE}$,
where each $y_\be\in L^2(0,1)$, equipped with the norm
$$
\|y\|^2_{L^2(\G)}=\sum_{\be\in\cE}\|y_\be\|^2_{L^2(0,1)}<\infty.
$$
We define the metric Laplacian $\D_M$ on $L^2(\G)$ by
$$
(\D_My)_\be=-y''_\be,\qqq y=(y_\be)_{\be\in\cE},
$$
where $(y''_\be)_{\be\in\cE}\in L^2(\G)$ and
$y$ satisfies the so-called Kirchhoff conditions:
\[
\lb{Dom1}
y \textrm{ is continuous on }\G,\qqq
\sum\limits_{\be=(v,\,u)_e\in\cE}\d_\be(v)\,y_\be'(v)=0, \qq \forall v\in V,
\]
$$
\d_\be(v)=\left\{
\begin{array}{rl}
1, \qq &\textrm{if $v$ is a terminal vertex of the edge $\be$, i.e. $t=1$ at $v$}, \\
 -1, \qq &\textrm{if $v$ is a initial vertex of the edge $\be$, i.e. $t=0$ at $v$}.
\end{array}\right.
$$

\subsection{Discrete Laplacians.}
Let $\ell^2(V)$ be the
Hilbert space of all square summable functions $f:V\to \C$, equipped
with the norm
$$
\|f\|^2_{\ell^2(V)}=\sum_{v\in V}|f(v)|^2<\infty.
$$
We define the self-adjoint normalized Laplacian (i.e., the Laplace operator) $\D$ on $\ell^2(V)$ by
\[
\lb{DOL}
 \big(\D f\big)(v)=
 -\frac1{\sqrt{\vk_v}}\sum\limits_{(v,\,u)_e\in\cE}\frac1{\sqrt{\vk_u}}\,f(u),
 \qquad v\in V,\qq  f\in \ell^2(V),
\]
where ${\vk}_v$ is the degree of the vertex $v\in V$ and all loops in the sum
\er{DOL} are counted twice.

We recall the basic facts about the Laplacian $\D$ (see
\cite{Ch97}, \cite{HS04}, \cite{MW89}), which
hold true for both finite and periodic graphs:

{\it (i) the point $-1$ belongs to the spectrum $\s(\D)$ and  $\s(\D)$ is contained in  $[-1,1]$, i.e.,
\[\lb{mp}
-1\in\s(\D)\subset[-1,1].
\]

(ii) On  periodic graphs the points
$\pm1$ are never eigenvalues of $\D$.}

We define \emph{the
fundamental graph} $\G_f=(V_f,\cE_f)$ of the periodic graph $\Gamma$
 as a graph
on the surface $\R^d/\Z^d$ by
\[
\lb{G0} \G_f=\G/{\Z}^d\ss\R^d/\Z^d.
\]
The fundamental graph $\G_f$ has the vertex set $V_f$ and the set
$\cE_f$  of unoriented edges, which are finite.
We identify the vertices of the fundamental graph $\G_f=(V_f,\cE_f)$
with the vertices of the graph $\G=(V,\cE)$ from the set $[0,1)^d$  by
\[
\lb{V0}
V_f=[0,1)^d\cap V=\{v_1,\ldots,v_\n\},\qqq \n=\# V_f<\infty,
\]
where $\n$ is the number of vertices of $\G_f$.
Denote by $\cB$ the set of all edges of $\G$ connecting the vertices
from $V_f$ with the vertices from $V\setminus V_f$. These edges will
be called \emph{bridges}.

The discrete Laplacian $\D$ on $\ell^2(V)$ has the standard decomposition into a
constant fiber direct integral  by
\[
\lb{raz}
\ell^2(V)={1\/(2\pi)^d}\int^\oplus_{\T^d}\ell^2(V_f)\,d\vt ,\qqq
U\D U^{-1}={1\/(2\pi)^d}\int^\oplus_{\T^d}\D(\vt)d\vt,
\]
$\T^d=\R^d/(2\pi\Z)^d$,
for some unitary operator $U$. Here $\ell^2(V_f)=\C^\nu$ is the fiber space
and $\D(\vt)$ is the Floquet  fiber $\nu\ts\nu$  matrix. The precise form of the matrix $\D(\vt)=\{\D_{jk}(\vt)\}_{j,k=1}^\n$
was determined in \cite{KS14}.

Each Floquet  matrix  $\D(\vt)$, $\vt\in\T^d$, has $\n$ real
eigenvalues $\l_n(\vt)$, $n\in\N_\n$, labeled by
\[
\label{eq.3}
\l_1(\vt)\leq\ldots\leq\l_{\nu}(\vt).
\]
 Each
$\l_n(\cdot)$, $n\in\N_\n$, is a continuous function on the
torus $\T^d$ and creates the spectral band $\s_n(\D)$ given by
\[
\lb{ban} \s_n=\s_n(\D)=[\l_n^-,\l_n^+]=\l_n(\T^d).
\]
Note that if
$\l_n(\cdot)= C_n=\const$ on some set $\mB\ss\T^d$ of positive
Lebesgue measure, then  the operator $\D$ on $\G$ has the eigenvalue
$C_n$ with infinite multiplicity. We call $C_n$ a \emph{flat band}.
The spectrum of the Laplace operator $\D$ on the periodic graph $\G$
has the form
\[
\lb{r0}
\s(\D)=\bigcup_{n=1}^\n\s_n(\D)=\s_{ac}(\D)\cup \s_{fb}(\D).
\]
Here $\s_{ac}(\D)$ is the absolutely continuous spectrum, which is a
union  of non-degenerated intervals from \er{ban}, and
$\s_{fb}(\D)=\{\m_1,\ldots,\m_r\}$, $r<\nu$,
is the set of
all flat bands (eigenvalues of infinite multiplicity). An open
interval between two neighboring non-degenerated bands is
called a \emph{gap}.
Let $\l_{\n_1}^+$, $\n_1\leq\n$, be the upper point of the absolutely continuous spectrum of the operator $\D$.  It is convenient for us an open
interval $(\l_{\n_1}^+,1)$
also to call a gap of the operator $\D$.

\subsection{Main results.}
Instead of the Laplacian $\D_M\ge 0$
it is convenient for us to define the momentum operator $\sqrt{\D_M}\ge 0$.
Due to Cattaneo Theorem (see Section 2) the spectrum of the operator $\sqrt{\D_M}\geq0$
on  $\G$ has the form
\[
\s(\sqrt{\D_M}\,)=\s_{ac}(\sqrt{\D_M}\,)\cup \s_{fb}(\sqrt{\D_M}\,).
\]
Both the sets
 $\s_{ac}(\sqrt{\D_M}\,)$ and $\s_{fb}(\sqrt{\D_M}\,)$
 are $2\pi$-periodic on the half-line $(0,\infty)$ and  are symmetric on the interval $(0,2\pi)$ with respect to the point $\pi$. Thus,
in order to study $\D_M$ it is sufficient to study its restriction $\O$ on the spectral interval $[0,\pi]$ given by
\[
\lb{O1} \O=\sqrt{\D_M}\,\chi_{[0,\pi]}(\sqrt{\D_M}\,)\;,
\]
where $\chi_A(\cdot)$ is the characteristic function of the set
$A$.
Due to Cattaneo Theorem (see Section~2) the spectrum of the operator $\O$ on a periodic metric graph $\G$ has  the form
\[\lb{Qr.1}
\begin{array}{c}
\displaystyle \s(\O)=\bigcup_{n=1}^\n\s_n(\O)=\s_{ac}(\O)\cup \s_{fb}(\O),\\[10pt]
\s_n(\O)=[z_n^-,z_n^+], \qqq -\cos(z_n^\pm)=\l_n^\pm,\qq n\in\N_{\n}.
\end{array}
\]
Here $\s_{ac}(\O)$ is a
union of non-degenerated spectral bands $\s_n(\O)$ with $z_n^-<z_n^+$ and $\s_{fb}(\O)$ is the flat band spectrum (for more details see Section 2).

\

%\subsection{Main results.}
Now we formulate our main results. Let
\[\lb{eq.b}
\textstyle \b=\sum\limits_{n=1}^\n\frac{\b_n}{\vk_n}\,,
\]
$\b_n$ is the bridge degree (the number of bridges incident to $v_n$)
and $\vk_n$ is the degree of $v_n\in \G_f$.

\begin{theorem}
\lb{TQL1}
i) All spectral bands $\s_n(\O)$ and $\s_n(\D)$, $n\in\N_\n$, of the momentum operator $\O$ and the discrete Laplacian $\D$, respectively, satisfy
\[\lb{sb1}
\textstyle|\s_n(\D)|\leq|\s_n(\O)|\leq{\pi\/\sqrt{2}}\,|\s_n(\D)|^{1\/2}.
\]

ii) The Lebesgue measure $|\s(\O)|$ and $|\s(\D)|$ of the spectrum of  \,$\O$ and $\D$, respectively, satisfies
\[
\lb{eq.7Q}
\textstyle|\s(\D)|\leq|\s(\O)|\leq{\pi\/\sqrt{2}}\,|\s(\D)|^{1\/2}\leq\pi\sqrt{\b}\;,
\]
where $\b$ is defined by \er{eq.b}.
Moreover, if there exist   spectral gaps $\g_1(\O),\ldots,\g_s(\O)$, in the spectrum $\s(\O)$,  then the following estimate holds true:
\[
\lb{GEga}
\begin{aligned}
\sum_{n=1}^s|\g_n(\O)|\ge \pi(1-\sqrt{\b}\,).
\end{aligned}
\]

iii) Let $\b<1$. Then the spectrum $\s(\D_M)$ has infinitely
many gaps.

\end{theorem}

\no \textbf{Remark.} 1) The estimate $|\s(\O)|\leq\pi\sqrt{\b}$ is not trivial iff
$\b<1$. The condition $\b<1$ holds true when the number of bridges in each vertex $v\in V_f$ is sufficiently small compared to the degree of the vertex. If we
change the coordinate system then the spectrum of the operator does not change, but the number of bridges does, in general. In order to get the best estimate in \er{eq.7Q} we have to choose a coordinate system when the number
$\b$ is minimal.

2) The Bethe-Sommerfeld conjecture states that each Schr\"odinger
operator $-\D+Q$ with a periodic potential $Q$ in $\R^d, d\geq2$ has only finitely many gaps in the spectrum. This conjecture was proved by Skriganov \cite{S85}. On an equilateral metric
graph  the spectrum of the Laplacian $\D_M$ has no gaps iff
$\s(\D)=[-1,1]$. If $\s(\D)\neq[-1,1]$, then in the spectrum of the
Laplacian $\D_M$ on a quantum graph there exist \textbf{infinitely
many gaps} $\g_1,\g_2,\ldots$ and $|\g_n|\rightarrow\infty$ as
$n\rightarrow\infty$.

\

{\bf Definition of Loop Graphs.} {\it i) A periodic graph $\G$ is
called  a loop  graph if each bridge $\be$ has the form
$\be=(v_j,v_j+\t(\be))$ for some $j\in\N_\n$ and $\t(\be)\in\Z^d$.

ii) A loop graph $\G$ is called a precise loop  graph if $\cos\lan\t
({\bf e}),\,\vt_0\ran=-1$ for all bridges $\be\in\cB$ and some
$\vt_0\in\T^d$. This point $\vt_0$ is called a precise point of the
loop  graph $\G$.}

\medskip

{\bf Remark.} 1) If  $\lan \vt_0,\t({\bf e})\ran/\pi $ is odd
for  all bridges $\be\in\cB$ and some vector $\vt_0\in\{0,\pi\}^d$,
then $\vt_0$ is a precise point of the loop graph $\G$.

2) The class of all precise loop graphs is large enough. The simplest
example  of precise loop graphs is the $d$-dimensional lattice. More
complicated examples are discussed in Proposition 2.3 in \cite{KS13}.

3) There exists a loop graph, which is not a precise loop graph. The simplest example of such graph is the triangular lattice (see Proposition 2.3 in \cite{KS13}).

\

We now describe all bands for precise loop periodic graphs.

\begin{theorem}\lb{QT100}
i) Let $\G$ be a loop graph. Then the spectral bands
$\s_n(\O)=[z_n^-,z_n^+]$  of the operator $\O$ satisfy
\[
\lb{Qeq.5} -\cos(z_n^-)=\l_n(0)=\l_n^-,\qqq \forall \;n\in \N_\n.
\]

ii) Let, in addition, $\G$ be  precise  with a precise point $\vt_0\in\T^d$. Then
\[
\lb{Qes}
 -\cos(z_n^+)=\l_n(\vt_0)=\l_n^+,\qqq \forall \;n\in \N_\n,
\]
\[
\lb{QQes.1}
\textstyle2\b=\sum\limits_{n=1}^\n|\s_n(\D)|
\leq\sum\limits_{n=1}^\n|\s_n(\O)|\,,\qqq %\b=\sum\limits_{n=1}^\n\frac{\b_n}{\vk_n}\,,
\]
where $\b$ is defined by \er{eq.b}.

\end{theorem}

We present the plan of our paper. In section \ref{Sec2} we estimate the
Lebesgue measure of the spectrum of the
operator $\sqrt{\D_M}$ on the finite interval $[0,\pi]$ in terms of geometric parameters of
the graph and discuss some spectral properties of metric Laplacians on loop graphs and bipartite periodic graphs. In Appendix we
collect spectral properties of the
discrete Laplacian from \cite{KS14}, needed to prove our main results.

%**********************************************
\section{\lb{Sec2} Proofs of the main theorems}
\setcounter{equation}{0}

\subsection{Cattaneo Correspondence.}

Cattaneo obtained a
correspondence between the spectrum of the Laplacian $\D_M$ on  the
equilateral metric graph and the spectrum of the Laplacian $\D$ on the
corresponding discrete graph \cite{C97}.
For the sake of completeness and the reader's convenience we recall this
correspondence.

Consider the eigenvalues problem with Dirichlet boundary conditions
\[
\lb{Dp}
-y''=E y,\qqq y(0)=y(1)=0.
\]
It is known that the spectrum of this problem  is given by $\s_D=\{(\pi n)^2 : n\in\N\}$. Here $(\pi n)^2$  is the so-called  Dirichlet eigenvalue of
the problem \er{Dp}.

We formulate Cattaneo's result \cite{C97} in the form convenient for us.
 This theorem gives a basis for describing
the spectrum  of the operator $\D_M$ in terms of $\D$, and conversely.

{\bf Theorem (Cattaneo)} \emph{i) The spectrum of the operator $\sqrt{\D_M}\geq0$
on the periodic metric graph $\G$ has the form
\[
\lb{al} \s(\sqrt{\D_M}\,)=\s_{ac}(\sqrt{\D_M}\,)\cup \s_{fb}(\sqrt{\D_M}\,),
\]
%where the absolutely continuous spectrum $\s_{ac}(\sqrt{\D_M}\,)$ and the flat %band spectrum $\s_{fb}(\sqrt{\D_M}\,)$ are given by
\[
\lb{MAb} \s_{ac}(\sqrt{\D_M}\,)=\big\{z\in\R_+\ : \ -\cos
z\in\s_{ac}(\D)\big\},
\]
\[
\lb{MFb} \s_{fb}(\sqrt{\D_M}\,)=\big\{z\in\R_+\ : \ -\cos
z\in\s_{fb}(\D)\big\}\cup\{\pi n : n\in\N\}.
\]}

\emph{ii) Each flat band $2\pi n$, $n\in\N$, is embedded in
the absolutely continuous spectrum $\s_{ac}(\sqrt{\D_M}\,)$.}

\emph{iii)  Both the sets
 $\s_{ac}(\sqrt{\D_M}\,)$ and $\s_{fb}(\sqrt{\D_M}\,)$
 are $2\pi$-periodic on the half-line $(0,\infty)$ and  are symmetric on the interval $(0,2\pi)$ with respect to the point $\pi$.}

\emph{iv) The spectrum of the operator $\O$ on a periodic metric
graph $\G$ has  the form
\[\lb{Qr}
\begin{array}{c}
\displaystyle \s(\O)=\bigcup_{n=1}^\n\s_n(\O)=\s_{ac}(\O)\cup \s_{fb}(\O),\\[10pt]
\s_n(\O)=[z_n^-,z_n^+], \qqq -\cos(z_n^\pm)=\l_n^\pm,\qq n\in\N_{\n}.
\end{array}
\]
Here $\s_{ac}(\D)$ is a
union of non-degenerated spectral bands $\s_n(\O)$ with $z_n^-<z_n^+$. Moreover, the first spectral band $\s_1(\O)=[0,z^+_1]$ is open. The flat band spectrum has the form
\[\lb{rel}
\s_{fb}(\O)=\{z_1,\ldots,z_r,\pi\},
\qqq
-\cos(z_k)=\m_k\neq1,\qqq k\in\N_r.
\]
}

\emph{v) $\s(\O)=[0,\pi]$ iff $\s(\D)=[-1,1]$.
}

\emph{vi) The spectrum of the operator $\O$ has exactly $k$ gaps iff $\s(\D)$
 has exactly $k$ gaps.}

\

\no \textbf{Remark.} 1) Cattaneo considered the Laplacian $\D_M$ on connected locally finite graphs (including periodic). In general, some points of the Dirichlet spectrum $\s_D$ are  not flat bands of the Laplacian $\D_M$.

2)  Von Below \cite{B85} considered the Laplacian $\D_M$ on connected finite graphs when the number of edges is grater than number of vertices.
For each $n\in\N$ he constructed an eigenfunction with the eigenvalue $(\pi n)^2$.
Without changing the proof this construction can be applied to a connected $\Z^d$-periodic graph ($d\geq2$) to obtain an eigenfunction with a compact support and the eigenvalue $(\pi n)^2$, $n\in\N$. Thus,
$\s_D\subset\s_{fb}(\D_M)$ (see \er{MFb}). Note that for $\Z$-periodic graph
some (all) points of the Dirichlet spectrum $\s_D$ may not be flat bands of the operator.

3) The relation between the spectra of $\D$ and $\sqrt{\D_M}$ is shown in
Fig.\ref{fRel}.

4) The flat bands $\pi n$,
$n\in\N$, of the operator $\sqrt{\D_M}\,$ will be called \emph{Dirichlet flat
bands}.

\setlength{\unitlength}{1.0mm}
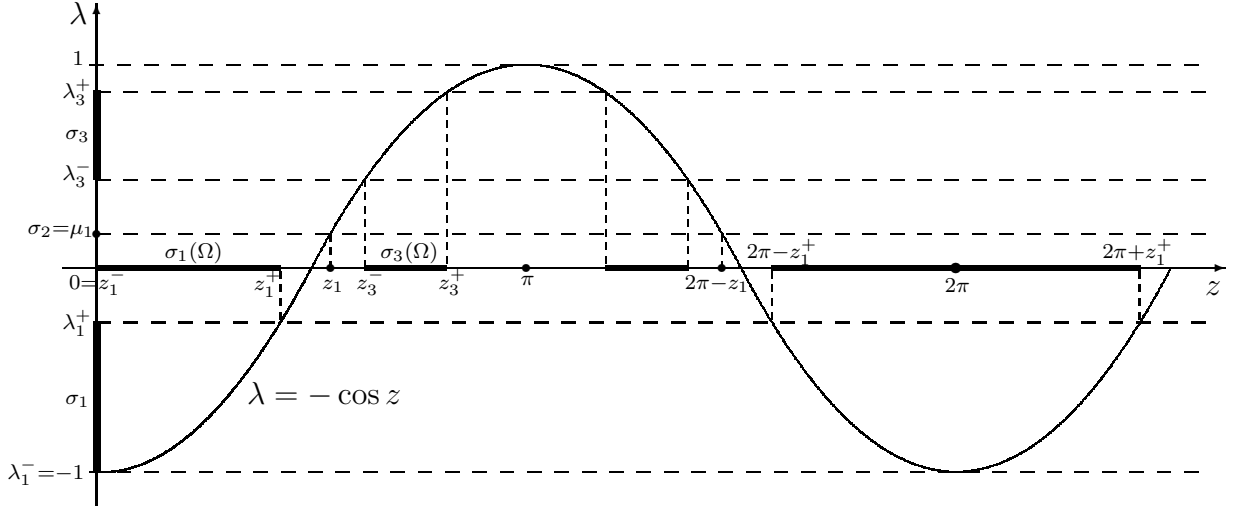
\begin{figure}[h]
\centering
\unitlength 0.9mm % = 2.845pt
\linethickness{0.4pt}
\ifx\plotpoint\undefined\newsavebox{\plotpoint}\fi % GNUPLOT compatibility
\begin{picture}(170,80)(0,0)

\put(-5,40){\vector(1,0){170.00}}

\put(0,5){\vector(0,1){74.00}}

\put(0.2,10){\line(0,1){22.20}}
\put(-0.2,10){\line(0,1){22.20}}
\put(0.4,10){\line(0,1){22.2}}
\put(-0.4,10){\line(0,1){22.20}}

\multiput(0,32)(4,0){41}{\line(1,0){2}}
\multiput(0,45)(4,0){41}{\line(1,0){2}}
\multiput(0,53)(4,0){41}{\line(1,0){2}}
\multiput(0,66)(4,0){41}{\line(1,0){2}}

\put(0,45){\circle*{1.0}}

\put(0.4,53){\line(0,1){13.20}}
\put(0.2,53){\line(0,1){13.20}}
\put(-0.2,53){\line(0,1){13.20}}
\put(-0.4,53){\line(0,1){13.20}}

\put(-13.0,9){$\scriptstyle\l_1^-=-1$}
\put(-5.0,31){$\scriptstyle\l_1^+$}
\put(-5.0,65){$\scriptstyle\l_3^+$}
\put(-5.0,53){$\scriptstyle\l_3^-$}
\put(-4.0,37.0){$\scriptstyle 0=z_1^-$}
\put(23.0,36.8){$\scriptstyle z_1^+$}
\put(95.0,41.5){$\scriptstyle 2\pi-z_1^+$}
\put(147.0,41.5){$\scriptstyle 2\pi+z_1^+$}
\put(10.0,41.5){$\scriptstyle \s_1(\Omega)$}
\put(41.5,41.5){$\scriptstyle \s_3(\Omega)$}
\put(22.0,20.0){$\l=-\cos z$}
\put(-10.0,45){$\scriptstyle\s_2=\m_1$}
\put(-4.5,59){$\scriptstyle\s_3$}
\put(-3.5,70){$\scriptstyle 1$}
\put(-4.5,20){$\scriptstyle\s_1$}
\put(-1,10){\line(1,0){2.00}}
\put(-1,70){\line(1,0){2.00}}

\multiput(-1,10)(4,0){41}{\line(1,0){2}}
\multiput(-1,70)(4,0){41}{\line(1,0){2}}
\bezier{600}(0,10)(15,9)(31.4,40)
\bezier{600}(94.2,40)(62.8,100)(31.4,40)
\bezier{600}(94.2,40)(125.6,-20)(157,40)

\put(34.2,40.0){\circle*{1.0}}
\put(33.0,37.0){$\scriptstyle z_1$}
\put(86.0,37.){$\scriptstyle 2\pi-z_1$}
%\put(75.5,41.5){$\scriptstyle 2\pi-\gS_3$}
\multiput(26.9,32.5)(0,2){4}{\line(0,1){1}}
\multiput(34.2,40.0)(0,2){3}{\line(0,1){1}}
\multiput(51.2,40.7)(0,2){13}{\line(0,1){1}}
\multiput(39.2,40.1)(0,2){7}{\line(0,1){1}}

\put(62.8,40.0){\circle*{1.0}}
\put(125.6,40.0){\circle*{1.5}}
\put(62.0,37.5){$\scriptstyle \pi$}
\put(124.0,36.5){$\scriptstyle 2\pi$}

\multiput(74.4,40.7)(0,2){13}{\line(0,1){1}}
\multiput(86.4,40.1)(0,2){7}{\line(0,1){1}}
\multiput(98.7,32.5)(0,2){4}{\line(0,1){1}}
\multiput(91.4,40.0)(0,2){3}{\line(0,1){1}}
\multiput(152.4,32.5)(0,2){4}{\line(0,1){1}}
\put(74.4,40.3){\line(1,0){12}}
\put(74.4,40.2){\line(1,0){12}}
\put(74.4,40.1){\line(1,0){12}}
\put(74.4,39.9){\line(1,0){12}}
\put(74.4,39.8){\line(1,0){12}}
\put(74.4,39.7){\line(1,0){12}}
\put(91.4,40.0){\circle*{1.0}}

\put(98.7,40.3){\line(1,0){53.8}}
\put(98.7,40.2){\line(1,0){53.8}}
\put(98.7,40.1){\line(1,0){53.8}}
\put(98.7,39.9){\line(1,0){53.8}}
\put(98.7,39.8){\line(1,0){53.8}}
\put(98.7,39.7){\line(1,0){53.8}}
\put(39.2,40.3){\line(1,0){12}}
\put(39.2,40.2){\line(1,0){12}}
\put(39.2,40.1){\line(1,0){12}}
\put(39.2,39.9){\line(1,0){12}}
\put(39.2,39.8){\line(1,0){12}}
\put(39.2,39.7){\line(1,0){12}}
\put(38.0,37.0){$\scriptstyle z_3^-$}
\put(50.0,37.0){$\scriptstyle z_3^+$}

\put(0,40.3){\line(1,0){26.9}}
\put(0,40.2){\line(1,0){26.9}}
\put(0,40.1){\line(1,0){26.9}}
\put(0,39.9){\line(1,0){26.9}}
\put(0,39.8){\line(1,0){26.9}}
\put(0,39.7){\line(1,0){26.9}}

\put(162.0,36.0){$z$}

\put(-4,76.0){$\l$}

\end{picture}
%\vspace{-0.5cm}
\caption{\footnotesize Relation between the spectra of $\D$ and $\sqrt{\D_M}$.}
\label{fRel}
\end{figure}

5) The number of flat bands of the operator $\O$
is $r+1$. Flat bands $z_1,\ldots,z_r$ corresponds to $r$ flat bands
of the discrete Laplacian and the flat band $\pi$ is a Dirichlet flat band.

6) Let $z_{\n_1}^+$, $\n_1\leq\n$, be the upper point of the absolutely continuous spectrum of the operator $\O$. If $z_{\n_1}^+<\pi$, then there is the gap $(z_{\n_1}^+,2\pi-z_{\n_1}^+)$ in the spectrum of $\sqrt{\D_M}$.
For the convenience the interval $(z_{\n_1}^+,\pi)$
is also called a gap of the operator $\O$.

\

%************************************************
\no {\bf Proof of Theorem \ref{TQL1}.} i) Consider the spectral band $\s_n(\D)=[\l_n^-,\l_n^+]$ of the discrete Laplacian $\D$ for some $n\in\N_\n$. Due to Cattaneo Theorem.iv the corresponding spectral band $\s_n(\O)$ of the momentum operator $\O$ has the form
\[\lb{f.1}
\s_n(\O)=[z_n^-,z_n^+], \qq \textrm{where} \qq \cos z_n^\pm=-\l_n^\pm.
\]
Applying Proposition \ref{mls}.i to the spectral band $\s_n(\D)$ and its preimage $\s_n(\O)$ under the function $\phi(z)=-\cos z$, $z\in[0,\pi]$, we obtain
\er{sb1}.

ii) The spectrum of the discrete Laplacian, defined by \er{r0}, is a union of non-overlapping segments $\gS_1,\ldots,\gS_{\n_\ast}$, where $\n_\ast\leq\n$, i.e.,
\[\lb{sp1}
\s(\D)=\bigcup\limits_{n=1}^{\n_\ast} \gS_n,\qqq
\gS_m\cap\gS_n=\varnothing, \qq m\neq n.
\]
Then, due to Cattaneo Theorem.iv, the spectrum of the operator $\O$ has the form
\[\lb{sp22}
\s(\O)=\bigcup\limits_{n=1}^{\n_\ast} \phi^{-1}(\gS_n),\qqq \textrm{where} \qq
\phi^{-1}(\gS_m)\cap \phi^{-1}(\gS_n)=\varnothing, \qq m\neq n.
\]
The identities \er{sp1}, \er{sp22} give
\[\lb{f.100}
|\s(\D)|=\sum_{n=1}^{\n_\ast}|\gS_n|,\qqq |\s(\O)|=\sum_{n=1}^{\n_\ast}|\phi^{-1}(\gS_n)|.
\]
Applying Proposition \ref{mls}.i to each segment $\gS_n$, $n\in\N_{\n_\ast}$, we obtain
\[
|\gS_n|\leq \big|\phi^{-1}(\gS_n)\big|,\qqq \forall n\in\N_{\n_\ast}.
\]
Summing these inequalities and using \er{f.100}, we obtain
\[
|\s(\D)|=\sum_{n=1}^{\n_\ast}|\gS_n|\leq\sum_{n=1}^{\n_\ast} \big|\phi^{-1}(\gS_n)\big|=|\s(\O)|,
\]
which gives the first inequality in \er{eq.7Q}.

We show the second inequality in \er{eq.7Q}.
Let $\gS_\ast\subset[-1,1]$ be the subset of the segment $[-1,1]$ with the Lebesgue measure $|\gS_\ast|=|\s(\D)|$, having the form
\[\lb{mL1}
\gS_\ast=[-1,-\l_\ast]\cup[\l_\ast,1],\qqq  \textstyle\l_\ast=1-{|\s(\D)|\/2}\,, \qq 0\leq\l_\ast<1.
\]
The preimage of the subset $\gS_\ast$ under the function $\phi(z)=-\cos z$, $z\in[0,\pi]$, is given by
\[
\phi^{-1}(\gS_\ast)=[0,z_\ast]\cup[\pi-z_\ast,\pi],\qqq \l_\ast=\cos z_\ast,\qqq \textstyle 0<z_\ast\leq{\pi\/2}\,.
\]
Combining Cattaneo Theorem.iv, Proposition \ref{mls} and the identity \er{mL1}, we obtain
$$
\textstyle|\s(\O)|=|\f^{-1}(\s(\D))|\leq|\f^{-1}(\gS_\ast)|
\leq{\pi\/\sqrt{2}}\,|\s(\D)|^{1\/2}.
$$
Thus, the second inequality in \er{eq.7Q} has been proved. The last inequality in \er{eq.7Q} follows from Theorem
\ref{T1}.i.

Now we will prove \er{GEga}.
Since $\s(\O)\subset[0,\pi]$, we have
\[\lb{GEga1}
\sum_{n=1}^s|\g_n(\O)|=
\pi-\big|\s(\O)\big|\geq \pi(1-\sqrt{\b}\,).
\]
Here we have used the estimate \er{eq.7Q}.

iii) Let $\b<1$. Then, due to \er{eq.7Q} we have
$|\s(\O)|<\pi$,
which yields that in the spectrum
of $\O$ there exists a gap. Hence, due to the periodicity of the spectrum of $\sqrt{\D_M}$, the spectrum
of $\D_M$ has infinitely
many gaps.

Now we give another proof of item iii). If $\b<1$, then Theorem \ref{T1}.i  gives that $|\s(\D)|<2$, i.e., $\s(\D)\neq[-1,1]$,
which, by Cattaneo Theorem.iii,v, yields that the spectrum
of $\D_M$ has infinitely
many gaps. \qq \BBox

\

\no {\bf Proof of Theorem \ref{QT100}.} i) Let $\G$ be a loop graph.
Then,  due to Theorem \ref{T100}.i, the spectral bands
$\s_n(\D)=[\l_n^-, \l_n^+]$ of the Laplacian $\D$ satisfy
$ \l_n^-=\l_n(0)$, $\forall\,n\in \N_\n$.
This and \er{Qr} give that the spectral bands $\s_n(\O)=[z_n^-,z_n^+]$
of  the operator $\O$ satisfy \er{Qeq.5}.

ii) Let $\G$ be a precise loop graph with a precise point $\vt_0\in\T^d$.
Cattaneo Theorem.iv and formula \er{es} give \er{Qes}.

Combining \er{es}, \er{es.1} and \er{Qes} with the simple estimate
$\cos z_n^--\cos z_n^+\le z_n^+-z_n^-$,
we obtain
$$
\textstyle2\b=\sum\limits_{n=1}^\n|\s_n(\D)|=
\sum\limits_{n=1}^\n(\l_n(\vt_0)-\l_n(0))=\sum\limits_{n=1}^\n(\cos z_n^--\cos z_n^+)\leq
\sum\limits_{n=1}^\n(z_n^+-z_n^-)=\sum\limits_{n=1}^\n\big|\s_n(\O)\big|.
$$
Thus, \er{QQes.1} has been proved. \qq
\BBox

\

A graph is called \emph{bipartite} if its vertex set is
divided into two disjoint sets (called \emph{parts} of the graph)
such that each edge connects vertices from distinct sets. It is known (\cite{BKS13}, \cite{HS04}) that \emph{a periodic graph is bipartite $\Leftrightarrow $ the point $1\in\s(\D)$ $\Leftrightarrow$ the  spectrum $\s(\D)$ is symmetric with respect to the point 0.}

Now we formulate some simple spectral properties of the Laplacian $\D_M$ and the momentum operator $\O$ on bipartite graphs.

\begin{theorem} \lb{TBG}
The following statements hold true.

i) If $\G$ is
bipartite,  then the flat bands $\pi^2(2n+1)^2$, $n\in\N$, of the Laplacian $\D_M$ are
embedded in $\s_{ac}(\D_M)$. If $\G$ is non-bipartite, then each
flat band $\pi^2(n+1)^2$, $n\in\N$, lies in a gap.

ii)  On the interval $(0,\pi)$ the spectrum of $\O$ is symmetric with respect
to the point $\pi/2$ iff $\G$ is bipartite iff the point $\pi\in\s_{ac}(\O)$.

iii) Let a fundamental graph $\G_f$ be bipartite. If $\nu$ is odd,
 then $z=\pi/2$ is a flat band of $\O$.

iv) Let $\G$ be a loop bipartite graph ($\G_f$ is not bipartite, since there
is  a loop on $\G_f$). Then each spectral band of the operator $\O$
on $\G$ has the form $\s_n(\O)=[z_n^{-},z_n^{+}]$, $n\in \N_\n$, where
\[
-\cos z^-_n=\l_n(0),\qqq \cos z^+_n=\l_{\n-n+1}(0).
\]

\end{theorem}

\no {\bf Proof.} i) If $\G$ is bipartite, then $1\in\s(\D)$. Since 1 is never a flat band, the point $1\in\s_{ac}(\D)$. Then \er{MAb} gives that
$\pi^2(n+1)^2\in\s_{ac}(\D_M)$ for all $n\in\N$. Similarly, if $\G$ is  non-bipartite, then $1\notin\s_{ac}(\D)$ and
$\pi^2(n+1)^2\notin\s_{ac}(\D_M)$ for all $n\in\N$. Thus, each
flat band $\pi^2(n+1)^2$, $n\in\N$, lies in a gap.

ii) Firstly, we show that the graph $\G$ is bipartite iff the point $\pi\in\s_{ac}(\O)$. Indeed, the graph $\G$ is bipartite iff the point $1\in\s_{ac}(\D)$. Due to \er{MAb}, the condition $1\in\s_{ac}(\D)$ is equivalent to $\pi\in\s_{ac}(\O)$.

Secondly, we prove that on the interval $(0,\pi)$ the spectrum of $\O$ is symmetric with respect to the point $\pi/2$ iff $\G$ is bipartite.
Let $\G$ be bipartite. Then the spectrum of the discrete Laplacian $\D$ is symmetric with respect to 0. Assume that $z\in(0,\pi)$ and $z\in\s(\O)$.
Then Cattaneo Theorem gives that $-\cos z\in\s(\D)$. Due to the symmetry of the spectrum $\s(\D)$, the point $-\cos(\pi-z)=\cos z\in\s(\D)$. This and \er{MAb} yield that $\pi-z\in\s(\O)$.

Conversely, let on the interval $(0,\pi)$ the spectrum of $\O$ be symmetric with respect to the point $\pi/2$. Since $-1\in\s_{ac}(\D)$, due to Cattaneo Theorem $0\in\s_{ac}(\O)$. Then the symmetry of $\s(\O)$ gives that $\pi\in\s_{ac}(\O)$ and $\G$ is bipartite.

iii) Let a fundamental graph $\G_f$ be bipartite and let $\nu$ be
odd.  By Theorem \ref{T1}.ii, $\m=0$ is a flat band of the discrete
Laplacian $\D$ on $\G$. This and Cattaneo Theorem.iv yield that
$z=\pi/2$ is a flat band of $\O$.

iv) By Theorem \ref{T100}.iii
each  spectral band of the discrete Laplacian on $\G$ has the form
$\s_n(\D)=[\l_n^{-},\l_n^{+}]$, $n\in \N_\n$, where $\l_n^{\pm}$ are
the eigenvalues of the matrix $\mp\D(0)$, i.e., $\l_n^-=\l_n(0)$,
$\l_n^+=-\l_{\n-n+1}(0)$. Then Cattaneo Theorem.iv gives that each
spectral band of the operator $\O$  has the form
$\s_n(\O)=[z_n^{-},z_n^{+}]$, $n\in \N_\n$, where $-\cos
z^-_n=\l_n^-=\l_n(0)$, $-\cos z^+_n=\l_n^+=-\l_{\n-n+1}(0)$. \qq \BBox

\

\no \textbf{Remark}. Item iii) gives a simple sufficient condition for existence of the flat band.

 %*******************************************
\section{\lb{Sec7}Appendix}
\setcounter{equation}{0}

We formulate some simple facts needed to prove our main results.

\begin{proposition}\lb{mls}
Let $\phi(z)=-\cos z$, $z\in[0,\pi]$. Then the following statements hold true.

i) If $\gS=[\l^-,\l^+]\subset[-1,1]$, then
\[\lb{emls}
\textstyle |\gS|\leq \big|\phi^{-1}(\gS)\big|\leq {\pi\/\sqrt{2}}\,|\gS|^{1\/2}.
\]

ii) If $\gS\subset[-1,1]$ is any subset of the segment $[-1,1]$ and  $\gS_\ast\subset[-1,1]$ is the subset with the same Lebesgue measure $|\gS_\ast|=|\gS|$, having the form
\[\lb{mL1.0}
\gS_\ast=[-1,-\l_\ast]\cup[\l_\ast,1],\qqq  \textstyle\l_\ast=1-{1\/2}\,|\gS|,
\]
then
\[\lb{Omm}
\textstyle|\phi^{-1}(\gS)|\leq|\phi^{-1}(\gS_\ast)|\leq{\pi\/\sqrt{2}}\,|\gS|^{1\/2}.
\]
\end{proposition}
\no {\bf Proof.} i) The function $\phi$ is an increasing bijection of the segment $[0,\pi]$ onto the segment $[-1,1]$. The preimage of the segment  $\gS=[\l^-,\l^+]\subset[-1,1]$ under the function $\phi$ has the form
\[\lb{f.1.0}
\phi^{-1}(\gS)=[z^-,z^+], \qq \textrm{where} \qq z^\pm=\phi^{-1}(\l^\pm).
\]
We have the simple inequality
\[\lb{f.4.1}
\begin{aligned}
\textstyle |\gS|=\l^+-\l^-=-\cos z^++\cos z^-=\int\limits_{z^-}^{z^+}\sin t\;dt \leq z^+-z^-=|\phi^{-1}(\gS)|.
\end{aligned}
\]
Thus, the first estimate in \er{emls} has been proved.

We will show that
\[\lb{f.2.0}
\big|\phi^{-1}(\gS)\big|\leq \big|\phi^{-1}(\gS_\ast)\big|,
\]
where $\gS_\ast\subset[-1,1]$ is given by
\[\lb{f.1.2}
\gS_\ast=[-1,-\l_\ast]\cup[\l_\ast,1],  \qqq \textstyle\l_\ast=1-{1\/2}\,(\l^+-\l^-)\in[0,1].
\]
Note that from the definitions of $\gS$ and $\gS_\ast$ it follows that
\[\lb{f.7.1}
|\gS_\ast|=|\gS|.
\]
The preimage of the segment  $\gS_\ast$ under the function $\phi$ has the form
\[\lb{f.1.1}
\textstyle \phi^{-1}(\gS_\ast)=[0,z_\ast]\cup[\pi-z_\ast,\pi], \qq \textrm{where} \qq z_\ast=\phi^{-1}(-\l_\ast)\in[0,{\pi\/2}].
\]
Using \er{f.1.0} and \er{f.1.1}
we can rewrite the inequality \er{f.2.0} in the equivalent form
\[\lb{f.2.00}
z^+-z^-\leq 2z_\ast.
\]
In order to prove \er{f.2.00} we use the definition of $z_\ast$ in \er{f.1.1}, \er{f.1.2} and consider the difference
$$
\begin{aligned}
\textstyle\cos {z^+-z^-\/2}-\cos z_\ast=\cos {z^+-z^-\/2}-
1+{1\/2}\,(\l^+-\l^-)=
\cos {z^+-z^-\/2}-1+{1\/2}\,(-\cos z^++\cos z^-)\\\textstyle=\cos{z^+\/2}\cos{z^-\/2}+\sin{z^+\/2}\sin{z^-\/2}
-\cos^2 {z^+\/2}-\sin^2{z^-\/2}\,.
\end{aligned}
$$
From this identity, using the simple inequalities
\[\lb{f.4.0}
\textstyle \cos{z^-\/2}\geq\cos{z^+\/2}\geq0,\qqq \sin{z^+\/2}\geq\sin{z^-\/2}\geq0, \qqq  0\leq z^-\leq  z^+\leq\pi,
\]
we obtain
\[\lb{f.5.0}
\textstyle\cos {z^+-z^-\/2}\ge\cos z_\ast \qq \Rightarrow \qq z_\ast\ge {z^+-z^-\/2}\,,
\]
since ${z^+-z^-\/2}\,,z_\ast\in[0,{\pi\/2}]$ and the function $\cos z$ decreases on this interval.
Thus, \er{f.2.00} (and, consequently, \er{f.2.0}) has been proved.
Using the estimate $\sin x\geq {2\sqrt{2}\/\pi}\,x$, $x\in[0,{\pi\/4}]$\,,
we obtain that the Lebesgue measure of the set $\gS_\ast$, defined by \er{f.1.2}, satisfies
\[
\lb{f.7.0}
\textstyle|\gS_\ast|=2\,(1-\l_\ast)=2\,(1-\cos z_\ast)=4\sin^2{z_\ast\/2}\geq{8\/\pi^2}\,z_\ast^2=
{2\/\pi^2}|\phi^{-1}(\gS_\ast)|^2,
\]
since ${z_\ast\/2}\in[0,{\pi\/4}]$.

Combining \er{f.2.0}, \er{f.7.0} and \er{f.7.1}, we obtain
\[\lb{f.8}
\big|\phi^{-1}(\gS)\big|^2\leq \big|\phi^{-1}(\gS_\ast)\big|^2=\textstyle 4z_\ast^2\leq
{\pi^2\/2}\,|\gS_\ast|={\pi^2\/2}\,|\gS|,
\]
which yields the second estimate in \er{emls}.

ii) The preimage $\f^{-1}(\gS_*)$ of the subset $\gS_\ast$ is given by
\er{f.1.1}.
Since the derivative $\f'$ of the function $\f(z)=-\cos z$ increases from the point $0$ to the point ${\pi\/2}$ and $\big|\phi^{-1}([-1,-\l_\ast])\big|=\big|\phi^{-1}([\l_\ast,1])\big|$, we obtain the first inequality in \er{Omm}.

Identity \er{f.7.1} and estimate \er{f.7.0}  imply  the second inequality in \er{Omm}.
\qq \BBox

\

We collect properties of the discrete Laplacian on periodic graphs
from \cite{KS14}, which we need to prove our main results.

\begin{theorem}\label{T1}

i) The Lebesgue measure
$|\s(\D)|$ of the spectrum of the Laplace operator $\D$ satisfies
\[
\lb{eq.7}
\textstyle |\s(\D)|\le \sum\limits_{n=1}^{\n}|\s_n(\D)|\le 2\b,
\]
where $\b$ is defined by \er{eq.b}.

ii) Let a fundamental graph $\G_f$ be bipartite. If $\nu$ is odd,
then $\m=0$ is a flat band of the discrete Laplacian $\D$ on $\G$.
\end{theorem}

\begin{theorem}\lb{T100}
i) Let $\G$ be a loop graph. Then the spectral bands $\s_n(\D)=[\l_n^-,
\l_n^+]$  of the Laplace operator $\D$ satisfy
\[
\lb{eq.5} \l_n^-=\l_n(0),\qqq \forall \;n\in \N_\n.
\]

ii) Let $\G$ be a precise loop graph with a precise point $\vt_0\in\T^d$.
Then
\[
\lb{es}
\s_n(\D)=[\l_n^-,\l_n^+]=[\l_n(0),\l_n(\vt_0)],\qqq \forall \;n\in \N_\n,
\]
\[
\lb{es.1}
\textstyle \sum\limits_{n=1}^\n|\s_n(\D)|=2\b,
\]
where $\b$ is defined by \er{eq.b}.

iii) Let in addition $\G$ be bipartite ($\G_f$ is not bipartite,
since  there is a loop on $\G_f$). Then each spectral band of the
Laplacian on $\G$ has the form $\s_n(\D)=[\l_n^{-},\l_n^{+}]$,
$n\in \N_\n$, where $\l_n^{\pm}$  are the eigenvalues of the matrix
$\mp\D(0)$.

\end{theorem}

%********************************************************

\medskip

\no\textbf{Acknowledgments.}
\footnotesize
Various parts of this paper were written during Evgeny Korotyaev's stay  in
the Mathematical Institute of Tsukuba University, Japan  and
Mittag-Leffler Institute, Sweden and Centre for Quantum Geometry of
Moduli spaces (QGM),  Aarhus University, Denmark. He is grateful to
the institutes for the hospitality. His study was partly supported
by The Ministry of education and science of Russian Federation,
project 07.09.2012 No 8501 and the RFFI grant "Spectral and
asymptotic methods for studying of the differential operators" No
11-01-00458 and the Danish National Research Foundation grant DNRF95
(Centre for Quantum Geometry of Moduli Spaces - QGM)".


\begin{thebibliography}{9999}
\setlength{\itemsep}{-\parskip}
\footnotesize

\bibitem[BKS13]{BKS13} Badanin A.; Korotyaev, E.; Saburova, N. Laplacians on $\Z^2$-periodic discrete graphs, preprint 2013.

\bibitem[B85]{B85} von Below, J. A characteristic equation associated to an eigenvalue problem on $c^2$-networks, Linear Algebra Appl. 71 (1985), 309-325.

\bibitem[BKu12]{BKu12} Berkolaiko G.; Kuchment, P. Introduction to Quantum
 Graphs, Mathematical Surveys and Monographs, V. 186 AMS, 2012.

\bibitem[C97]{C97} Cattaneo, C. The spectrum of the continuous Laplacian
 on a graph, Monatsh. Math. 124 (1997), 215--235.

\bibitem[Ch97]{Ch97} Chung, F. Spectral graph theory, AMS,
 Providence, Rhode. Island, 1997.

\bibitem[HS04]{HS04} Higuchi, Y.; Shirai, T. Some spectral and
 geometric properties for infinite graphs, AMS Contemp. Math. 347 (2004), 29--56.

\bibitem[K98]{K98}  Korotyaev, E. Estimates of periodic potentials in terms
of gap lengths. Comm. Math. Phys. 197 (1998), no. 3, 521--526.

\bibitem[K03]{K03} Korotyaev, E. Characterization of the spectrum of
 Schr\"odinger operators with periodic distributions.
 Int. Math. Res. Not. 2003, no. 37, 2019--2031.

\bibitem[KS13]{KS13} Korotyaev, E.; Saburova, N. Schr\"odinger operators on
 periodic discrete graphs, preprint 2013.

\bibitem[KS14]{KS14} Korotyaev, E.; Saburova, N. Spectral estimates for
normalized Laplacian and its perturbations on periodic discrete
graphs, preprint 2014.

\bibitem[LP08]{LP08} Lled\'o, F.; Post, O. Eigenvalue bracketing for
 discrete and metric graphs, J. Math. Anal. Appl. 348 (2008), 806--833.

\bibitem[MW89]{MW89} Mohar, B.; Woess, W. A survey on spectra of infinite graphs,
Bull. London Math. Soc., 21 (1989), 209--234.

\bibitem[P12]{P12} Post, O. Spectral analysis on graph-like spaces.
 Lecture Notes in Mathematics, 2039. Springer, Heidelberg, 2012.

\bibitem[S85]{S85}  Skriganov, M. The spectrum band structure of
 three-dimensional Schr{\"o}dinger operator with periodic potential,
  Invent. Math. 80 (1985), 107--121.


\end{thebibliography}
\end{document}